\documentclass[a4paper,12pt]{article}
\usepackage{epsf,amsfonts,diss,mydefs,url}%,hyperref,showkeys}
\def\g{\gamma}

\def\paper{paper}
\def\thep{.}
\newcommand{\one}{\relax\ifmmode{\rm 1\>\!\!\!I}\else{$\rm 1\!I$}\fi}
\def\computation#1{#1}
\newcommand{\crt}{\mathop{\rm cr}}
\begin{document}
%%%
%%% Some more defs:
%%%
\def\q{Q}

%\title{On discrete curves in $\C^2$ and $\CP^1$.}
\title{Discrete curves in $\CP^1$ and the Toda lattice.}
\author{
  Tim Hoffmann\footnote{Supported by the DFG sfb288. 
    email: \protect\url{hoffmann@math.tu-berlin.de}}
 \  and Nadja Kutz\footnote{Supported by the Servicezentrum Humboldt
   University of Berlin.
    email: \protect\url{nadja@math.tu-berlin.de}}
  }
\maketitle
\label{sec:dc}
\begin{abstract}
  In this \paper\ we investigate flows on discrete curves
  in $\C^2$, $\CP^1$, and $\C$. A novel interpretation of 
  the one dimensional Toda lattice hierarchy and reductions thereof
  as flows on discrete curves will be given.
\end{abstract}
\section{Introduction}
In this \paper\ we investigate flows on discrete curves
in $\C^2$, $\CP^1$, and $\C$. 
By a discrete curve $\g$ in $K$ we mean a map $\g:\Z \rightarrow K$.
A flow on $\g$ is a smooth variation of $\g$.
The description of these curves is mainly motivated by the
picture of discrete curves in $\CP^1$, in other words we 
imagine curves in $\C^2$ as curves in $\CP^1$, which 
are lifted to homogeneous coordinates; and curves in $\C$, as curves in
$\CP^1$, which do not hit infinity. 
In particular this view allows us to find a novel interpretation of 
the wellknown one dimensional Toda lattice hierarchy in terms of
flows on discrete curves. The Toda lattice hierarchy is a set of
equations, including the Toda equation:
$$
\ddot q_k = e^{q_{k+1}-q_k}-e^{q_{k}-q_{k-1}}.
$$
The Toda equation is sometimes also called first flow equation of the
Toda lattice hierarchy, it was discovered by Toda in 1967 (\cite{TO89}). 
A good overview about the vast literature about the Toda lattice
can be found in (\cite{FT86,SU02}).
The paper is organized as follows. 

In section \ref{sec:dc:definitions} we define discrete curves in $\C^2$
and $\CP^1.$ A discrete analog of the Schwartzian derivative in
terms of cross-ratios of four neighboring points will be defined. 
A zero curvature or Lax representation for
flows on that discrete curves will be given.

In section \ref{sec:dc:arclength} we will restrict ourselves to the case
of so-called discrete (conformal) arc length parametrized curves. 
In \ref{subsec:dc:arclengthC2} we define certain flows on
 arc length parametrized curves which in the turn induce flows for the
cross-ratios of that curve. It will be shown that for  these flows  the 
cross-ratios give solutions to equations of the
Volterra hierarchy.

In \ref{subsec:dc:discreteEuclidean} the arc length parametrized curves in
$\C^2$ viewed as curves in $\CP^1$ in homogeneous coordinates will be
reduced to curves in $\C$, by assuming that the second coordinate does
not vanish. We call this reduction Euclidean reduction.  It will be
shown that if one starts with a flow on a discrete arc length
parametrized curve in $\C^2$, such that the corresponding cross-ratios
are solutions to the second flow of the Volterra hierarchy, which is
also called the discrete KdV flow, then the corresponding curve
obtained by Euclidean reduction satisfies a discrete version of the
mKdV flow. This gives a discrete version of a geometrical
interpretation of the Miura transformation.

In \ref{subsec:dc:discreteFlows} a nice geometrical interpretation
of B\"acklund transformations of arc length parametrized discrete
curves  in  $\C^2$ will be given. In the turn a time discrete
version of the Volterra equation, which is the first flow equation
of the Volterra hierarchy, is derived.

In section~\ref{sec:dc:todaLattice} determinants of two and three neighboring
points of a discrete curve in $\C^2$ will be identified with the
Flaschka-Manakov variables of the Toda lattice. The Toda flows
define flow directions for curves in $\C^2.$ 

In \ref{sec:dc:threeFlowsRed} the flows on the discrete curves in $\C^2$ 
given by  Toda lattice hierarchy will be further investigated.
In particular in \ref{subsec:dc:threeFlows} the flows corresponding
to the ''first'' three flows will be given explicitly. 

In  \ref{sec:dc:invariant} we derive several
geometrical features of these flows. In particular in
\ref{subsec:dc:ellipse} we will
look at flows on curves which are compatible with the
 known reduction   \cite{SU02} of the Toda lattice hierarchy
to the Volterra lattice hierarchy by setting the Flaschka-Manakov  variable
$p_k= 0.$ This reduction is different from from the reduction
in \ref{subsec:dc:arclengthC2}, where among others $g_k = {\exp(({q_k-q_{k+1}})/{2})} = const$.
Nevertheless the cross ratios
of such curves evolve again with equations of the Volterra hierachy. Hence
there exist two geometrically motivated reductions of the Toda lattice
hierachy to the Volterra hierarchy.
The case $p_k= 0$ and $g_k= const$ lies in the 
intersection of both reductions, it belongs to the trivial
solutions of the Toda lattice hierachy
equations. The corresponding class of curves, which are
invariant under this trivial flow are quadrics, as will be shown
in that subsection. 

In the subsection~\ref{sec:dc:elastica} we will again look for
a class of curves, which is invariant (up to euclidian motion
and tangential flow)  under a certain Toda flow. This time
we find that in the Euclidian reduction socalled discrete
planar elastic curves are invariant under the discrete mKdV flow.

%%%
%%% Flows on discrete curves in CP1
%%%
\section{Flows on discrete curves} \label{sec:dc:definitions} 
%in $\protect\CP^1$}
%We restrict ourselves to planar discrete curves for a while.
\subsection{Flows on discrete curves in $\C^2$}
Let $c:\Z\to\CP^1$, $k\mapsto c_k$ be a discrete curve in the complex projective space.
We assume $c$ is immersed, i.e. $c_{k-1},c_k$ and $c_{k+1}$ are pairwise disjoined. 
By introducing homogeneous coordinates, we can lift $c$ to a map 

\beqna \g:\Z&\to&\C^2\nonumber\\
k&\mapsto& \g_k = \left(\ba{c} x_k \\ y_k \ea\right)
\label{def:dc:gamdef}\eeqna with $c_k= x_k y_k^{-1}$. Obviously $\g$
is not uniquely defined: For $\lambda:\Z\to\C^*$, $\lambda_k\g_k$ is also
a valid lift. 

\noindent Define:
\beqna
g_k &=& \det(\g_k,\g_{k+1}) = x_k y_{k+1}-y_k x_{k+1}\\
u_k&=&\det(\g_{k-1},\g_{k+1})  \label{def:dc:gu}
\eeqna
The following lemma can be straightforwardly obtained by using
the above definitions (\ref{def:dc:gu}):
\begin{lemma}
$$\g_{k+1}=\frac{1}{g_{k-1}}(u_k \g_k - g_k \g_{k-1}).
\label{lem:dc:Gam1}
$$
\label{lem:dc:CurveSteps}\end{lemma}
If the  variables $u_k$ and $g_k$ and initial points
$\g_0$ and $\g_1$ are given, then lemma~\ref{lem:dc:Gam1} is
a recursive definition of a discrete curve. This may look
as an odd way to define a discrete curve;
nevertheless later on the power of this definition
will become more apparent. In particular the variables $u_k$ and $g_k$
will be related to the Flaschka-Manakov \cite{Flasch74a,Flasch74b,MAN74} 
variables  of the
one dimensional Toda lattice.

Note that after choosing  an initial $\g_0$
it is always possible (c is immersed) to find $\g_k$ such that:
\beqn
g_k = 1 \label{def:dc:detEins}
\eeqn
for all $k$.
We will call discrete curves $ \g$ with property (\ref{def:dc:detEins})  
conformal arc length parametrized.
Hence the variables $g_k$ measure the deviation from 
arc length parameterization.      
\begin{definition}
  The {\em cross-ratio}\/ of four points $a,b,c,d\in \C^2$ is defined by
  \[\crt(a,b,c,d) =
  \frac{\det(a,b)}{\det(b,c)}\frac{\det(c,d)}{\det(d,a)}\thep\]
\end{definition}
Let us denote the cross-ratio of four neighboring points of $\g$ by $\q$:
\begin{equation}
  \label{eq:dc:q}
  \q_k :=
  \crt(\g_{k-1},\g_k,\g_{k+2},\g_{k+1})=\frac{g_{k-1}g_{k+1}}{u_k u_{k+1}}\thep
\end{equation}

Now $\det(\g_k, \frac{\g_{k+1}-\g_{k-1}}{g_k + g_{k-1}})= 1$. This
means that  
$\g_k$ and  $\frac{\g_{k+1}-\g_{k-1}}{g_k + g_{k-1}}$
 are linearly independent. 
Hence an arbitrary flow on (or variation of) $\g$  can be written in the 
following way:
\begin{equation}
       \frac{d}{dt}\g_k= \dot\g_k = \alpha_k \g_k + \frac{\beta_k}{u_k}(\g_{k+1}-\g_{k-1})
       \hspace{4mm}\alpha_k,\beta_k \in \C \thep\label{eq:dc:generalFlow}
\end{equation}
The variables $\alpha_k,\beta_k \in \C $ are arbitrary, for
convenience we sometimes use the normalized variable:
\beqn
\hat \beta_k = \beta_k (g_{k+1}+g_{k-1}).
\label{def:dc:betahut}
\eeqn

\begin{lemma}
 A flow on the discrete curve $\g$ given by (\ref{eq:dc:generalFlow}) 
 generates the following flow on the variables $g_k$:
\begin{equation}
       \dot g_k = g_k (\alpha_{k+1}+\alpha_k) +\beta_{k+1}-\beta_k\thep
        \label{eq:dc:dotg}
\end{equation}
\label{lem:dc:gu}\end{lemma}
\begin{proof}
  Differentiate $g_k$ in (\ref{def:dc:gu}): 
  \begin{equation}
    \begin{array}{rcl}
      \dot g_k &=&\det(\dot\g_k,\g_{k+1}) + \det(\g_k,\dot\g_{k+1})  \\
      \computation{%
        &=&\det(\alpha_k\g_k+\frac{\beta_k}{u_k}(\g_{k+1}-\g_{k-1}),\g_{k+1})\\ 
 &&+ 
        \det(\g_k,\alpha_{k+1}\g_{k+1} + \frac{\beta_{k+1}}{u_{k+1}}(\g_{k+2}-\g_k))\\
        }%
      &=& g_k \alpha - \beta_{k} + g_k \alpha_{k+1} + \beta_{k+1}\thep
    \end{array}
%    \label{}
  \end{equation}

\end{proof}
If the flow is arc length preserving, i.e. in particular $\dot g_k=0$
f.a. $k \in \Z$ then there exists an obvious trivial
solution for (\ref{eq:dc:dotg}): Choosing $\beta\equiv 0$ induces 
$\alpha_{k+1} = -\alpha_k$. This flow corresponds to the freedom of the initial
choice of $\g_0$ and has no effect on the corresponding curve $c$
in $\CP^1$. Note also that in this case
(\ref{eq:dc:dotg}) is a linear equation. So one can
always add any two flows solving it.

\begin{lemma}
A flow on the discrete curve  $\g$ given by (\ref{eq:dc:generalFlow}) 
 generates the following flow on the variables $u_k$:
\beqn\begin{array}{rcl}
\dot u_k &=& u_k(\alpha_{k-1}+ \alpha_{k+1})\\&+&
 \beta_{k-1}\frac{g_{k}}{u_{k-1} g_{k-1}}(g_{k-2}+g_{k-1}) -
\beta_{k+1}\frac{g_{k-1}}{u_{k+1} g_{k}}(g_{k}+g_{k+1})\\
&+&u_k (\frac{1}{g_{k}}\beta_{k+1} - \frac{1}{g_{k-1}}\beta_{k-1}) 
\label{eq:dc:FlowOnU}
\end{array}
\eeqn
or equivalently
\beqn\begin{array}{rcl}
\frac{\dot u_k}{u_k} &=& \alpha_{k-1}+ \alpha_{k+1}\\&&+
 \beta_{k-1}\frac{\q_{k-1}}{g_{k-2} g_{k-1}}(g_{k-2}+g_{k-1}) -
\beta_{k+1}\frac{\q_{k}}{g_{k} g_{k+1}}(g_{k}+g_{k+1})\\
&&+\frac{1}{g_{k}}\beta_{k+1} - \frac{1}{g_{k-1}}\beta_{k-1}
\end{array}
\eeqn
 \label{lem:dc:flowOnU}\end{lemma}
\begin{proof}
\beqn
\ba{rcl}\dot u_k  &=& \det(\dot \g_{k-1}, \g_{k+1})+ 
\det(\g_{k-1},\dot \g_{k+1})\\
&=& \alpha_{k-1}u_k + g_k \frac{\beta_{k-1}}{u_{k-1}}-
\frac{\beta_{k-1}}{u_{k-1}}\det(\g_{k-2},\g_{k+1})\\
 &&+ \alpha_{k+1}u_k - g_{k-1} \frac{\beta_{k+1}}{u_{k+1}}+
\frac{\beta_{k+1}}{u_{k+1}}\det(\g_{k-1},\g_{k+2})
\ea \eeqn
Using lemma \ref{lem:dc:CurveSteps} one gets
$$
\det(\g_{k-1},\g_{k+2})= \frac{1}{g_k}u_k u_{k+1}
- \frac{g_{k+1} g_{k-1}}{g_k}$$ which gives the result.\end{proof}

\begin{lemma}
A flow on the discrete curve  $\g$ given by (\ref{eq:dc:generalFlow}) 
 generates the following flow on the cross-ratios $\q_k$:
\beqn
\ba{rcl}
\frac{\dot \q_k}{\q_k}&=&(\q_k
 -1)(\beta_{k+1}(\frac{1}{g_{k+1}}+\frac{1}{g_{k}})
-\beta_k (\frac{1}{g_{k}}+\frac{1}{g_{k-1}}))\\
&&+\q_{k+1} \beta_{k+2}(\frac{1}{g_{k+2}}+\frac{1}{g_{k+1}})\\
&&-\q_{k-1} \beta_{k-1}(\frac{1}{g_{k-1}}+\frac{1}{g_{k-2}})
\ea
\label{eq:dc:FlowOnQ}\eeqn
\end{lemma}
\begin{proof}
Using the definition of the cross-ratio (\ref{eq:dc:q}) 
and lemmas \ref{lem:dc:flowOnU} and \ref{lem:dc:gu} the assertion
follows immediately.
\end{proof}
\subsection{Flows on discrete curves in $\CP^1$}
\begin{lemma}\label{thm:dc:cp1Flow}
  A flow on the discrete curve  $\g$ given by (\ref{eq:dc:generalFlow}) 
  generates the following flow on the non-lifted curve $c$ in $\CP^1$
  whenever $c$ does not hit $\infty$:
  \begin{equation}
    \label{eq:dc:flowOnCP1}
    \dot c_k = \beta_k \frac{g_k + g_{k-1}}{2g_k g_{k-1}}\; 2
    \frac{(c_{k+1}-c_k)(c_k-c_{k-1})}{c_{k+1}-c_{k-1}}\thep
  \end{equation}
\end{lemma}
Equation~(\ref{eq:dc:flowOnCP1}) can be written in the form 
\[
\dot c_k = \frac{\beta_k}{M^h(g_k,g_{k+1})} M^h(c_{k+1}-c_k,c_k-c_{k-1})
\]
where $M^h$ denotes the harmonic mean. In general the flow of the
non-lifted curve will depend on the chosen lift, since the $g_k$
depend on the choice. Note however, that the
$\alpha_k$ do not contribute to the evolution of the non-lifted curve
$c$. 

\noindent\begin{proof}
  For a given lift $\g={x_k \choose y_k}$ the non-lifted curve may be
  reconstructed by $c_k = x_k/y_k$ whenever $y_k\neq0$ which means
  that the curve does not hit $\infty$. Now insert this and the
  definitions of $g$ and $u$ in both sides of
  equation~(\ref{eq:dc:flowOnCP1}).
\end{proof}
\subsection{Zero Curvature representation}
Define:
\beqn
F_k
=\left(\ba{l}\g^T_{k}\\\g^T_{k-1}\ea\right)=
\left(\ba{ll}x_k&y_k\\x_{k-1}&y_{k-1}\ea \right).
\eeqn
Note that if the curve $\g$ is conformal arc length parametrized
(\ref{def:dc:detEins}) then $F_k \in SL(2,\C)$ for all $k \in \Z$.
\begin{Proposition}
Let $\alpha_k, \beta_k \in \C$ be arbitrary and $g_k$,$u_k$ be as
defined in (\ref{def:dc:gu}). Then 
$$F_{k+1}=L_k F_k \hspace{5mm} \dot F_k = V_k F_k$$
with
\beqn
L_k = \left(
  \ba{cc}\frac{1}{g_{k-1}}u_k&-\frac{g_k}{g_{k-1}}\\1&0\ea\right)
\hspace{5mm}
V_k = \left(\ba{cc}
\alpha_{k}+\frac{1}{g_{k-1}}\beta_{k}&
-(1+\frac{g_{k}}{g_{k-1}})\frac{\beta_{k}}{u_{k}}\\
(1+\frac{g_{k-2}}{g_{k-1}})\frac{\beta_{k-1}}{u_{k-1}}&
\alpha_{k-1}-\frac{1}{g_{k-1}}\beta_{k-1}
\ea \right).
\label{eq:dc:vk}\eeqn
The compatibility equation
\beqn
\dot L_k = V_{k+1}L_k - L_k V_k \label{eq:dc:comp}
\eeqn
is satisfied for all $\alpha_k, \beta_k \in \C$. 
\label{prop:dc:LV}
\end{Proposition}
The compatibility equation \ref{eq:dc:comp} is also called
zero curvature equation.

\noindent\begin{proof}
The construction of $L_k$ is obvious with lemma
\ref{lem:dc:CurveSteps}.
The construction of $V_k$ follows also quite
straightforwardly from lemma
\ref{lem:dc:CurveSteps} and the definition of
the flow on $\g$ in (\ref{eq:dc:generalFlow}).
The compatibility equation (\ref{eq:dc:comp}) holds 
by construction. The flows on $g_k$ and $u_k$
were constructed by using a well defined flow on $\g$
for which in particular $\frac{d}{dt}( L_k \g_{k})= V_{k+1}\g_{k+1}$
(which gives (\ref{eq:dc:comp})).
Nevertheless (\ref{eq:dc:comp})
can also easily be checked directly.
\end{proof}

%%%
%%% The case g == 1
%%%
%\section{Discrete curves in complex projective space}
\section{Conformal arc length parametrized curves } %\protect{$\C^2$}}
\label{sec:dc:arclength}
\subsection{Conformal arc length parametrized curves in
  \protect{$\C^2$} }
\label{subsec:dc:arclengthC2}
Let $c$ be a discrete curve in $\CP^1$. Up M\"obius transformations
(or up to the choice of $c_0$, $c_1$, and $c_2$)
$c$ is completely determined by the cross-ratios $\q_k$. If one scales
all $\q_k$ with a non-vanishing factor $\lambda$ one gets---again up to
M\"obius transformations---a new
discrete curve $c(\lambda)$. We call the family of all such curves the
{\em associated family} of $c = c(1)$.

As mentioned in the beginning, one has a choice when lifting a
discrete curve from $\CP^1$ to $\C^2$. This choice can be fixed by
prescribing the determinants of successive points. A natural coice is
here to set $g\equiv 1$ which we called conformal arc length
parametrization and we will discuss this choice here. However, in 
section~\ref{sec:dc:todaLattice} we will see, that other
normalizations are likewise meaningfull.

%In this section we assume, that all $g_k=1$. 
We will now investigate flows on $\g$ that preserve the conformal arc
length. The condition for this is, that $\dot g_k = 0$, which
implies by equation~\ref{eq:dc:dotg}
\begin{equation}
  \label{eq:dc:arcLengthConstraint}
  \alpha_{k+1}+\alpha_k = \beta_{k+1}- \beta_k
\end{equation}
for $\alpha$ and $\beta$ from equation~(\ref{eq:dc:generalFlow}).

Recall that the conformal arc length condition leaves us with an
initial choice of $\g_k$ and that this freedom corresponds to
a trivial flow with $\beta\equiv0$ and $\alpha_{k+1} =
-\alpha_k$. This flow is a first example of a (conformal) arc length
preserving flow. It does not change the curve $c$ in $CP^1$ though,
since only the $\beta_k$ contribute to the evolution of the non-lifted
curve.

If we choose $\beta\equiv 1/2$ and $\alpha \equiv 0$. We get for the curve
\begin{equation}
        \dot\g_k = \frac1{2 u_k}(\g_{k+1}-\g_{k-1})\thep
        \label{eq:dc:conformalTangentialFlow}
\end{equation}
This is what we will call the {\em conformal tangential flow}.  Then
$\dot u_k = 
\frac1{u_{k-1}} -\frac1{u_{k+1}}$ and $\q$ will solve the famous
Volterra model \cite{FT86,SU02}:
\begin{equation}
  \label{eq:dc:Volterra}
  \dot \q_k = \q_k(\q_{k+1} - \q_{k-1})\thep
\end{equation}
If we want this equation for the whole associated family of $\g$ we
must scale time by $\lambda$:
\[\lambda\dot \q_k(\lambda) = \q_k(\lambda)(\q_{k+1}(\lambda) -
\q_{k-1}(\lambda))\] 
One obtains the next higher flow of the Volterra hierarchy \cite{SU02}
when one chooses 
%$\beta =1/2  2\M \q +1 = 1/2(\q_{k+1}+\q_{k}+1)$. 
$\beta_k = 1/2(\q_{k-1}+\q_{k}+1)$. 
This implies
\begin{equation}
  \label{eq:dc:Volterra2}
  \dot \q_k = \q_k(\q_{k+1}(\q_{k+2} + \q_{k+1} + \q_k) - \q_{k-1}
  (\q_k + \q_{k-1} + \q_{k-2}))\thep 
\end{equation}

\noindent {\bf Conjecture}

{\em
\noindent Let
\beqn
\beta^{\mbox{new}}_{k+1}-\beta^{\mbox{new}}_{k} := 
\frac{\dot \q^{\mbox{old}}_k}{\q^{\mbox{old}}_k}.
\label{def:dc:betaIt}\eeqn
Given flows $\frac{\dot \q_k}{\q_k}$ this defines the variables 
$\beta^{\mbox{new}}_{k}$ up to a constant. Now observe that starting with the
flow $\frac{\dot \q_k}{\q_k}=0$ gives $\beta_k=a_1$, where $a_1 \in
\C$
is an arbitrary constant. Inserting these $\beta_k$ into the
flow equation (\ref{eq:dc:FlowOnQ}) in the reduced case $g_k=1$
gives in the turn a new flow equation
$$
\frac{\dot \q_k}{\q_k}=a_1(\q_{k+1}-\q_{k-1}),
$$
which is (up to the constant $a_1$) the Volterra equation.
Now inserting this Volterra equation into 
equation (\ref{def:dc:betaIt}) gives new $\beta_k$ as
$$
\beta_k := a_1(\q_k + \q_{k-1})+ a_2
$$
which in the turn give the following flows on the cross-ratios:
\beqnao
\frac{\dot \q_k}{\q_k}
&=&a_1 \left( (\q_k-1)(\q_{k+1}-\q_{k-1})\right)\\
&+&a_1 \left(\q_{k+1}(\q_{k+2}+\q_{k+1})-\q_{k-1}(\q_{k-1}+\q_{k-2})\right)\\ 
&+& a_2(\q_{k+1}-\q_{k-1}).
\eeqnao
This is (up the constant $a_1$) the next higher flow in the
Volterra hierarchy plus a Volterra term if $a_2 \neq 0$. 

We conjecture that all higher flows of the Volterra hierarchy can be
obtained in this way. There are strong indications that this
holds also in the continuous case \cite{UP}. A proof of this 
conjecture would be beyond the scope of this article, we postpone
this to a later publication.
}

\noindent
To make contact with the classical results we will now derive the
$2\times 2$-Lax representation of the Volterra model for our
tangential flow:

Define the gauge matrix 
\begin{equation}
  \label{eq:dc:volerraGaugeE}
  E_k := \prod_{i=0}^{k-1} u_i \quadmatrix{1}{0}{0}{u_k}
\end{equation}
and set 
\begin{equation}
  \label{eq:dc:VolterraL}
  \tilde L_k := E_{k+1}^{-1}L_kE_k = \quadmatrix{1}{-1}{\q_k}{0}
\end{equation}
This gauge of $L_k$ implies the following change for $V_k$:
\begin{equation}
  \label{eq:dc:volterraV}
  \begin{array}{rcl}
    \tilde V_k &:=& E_k^{-1}V_k E_k + E_k^{-1}\dot E_k\\[0.2cm]
 &=& \quadmatrix{1+\q_{k-1}}{-1}{\q_{k-1}}{\q_k} + \quadmatrix{
\q_{-1}-\frac12}{0}{0}{\q_{-1}-\frac12}\thep
\end{array}
\end{equation}
The second matrix summand may be omitted since it is constant.
If we now transpose the system, reverse the direction of the
$k$-labeling and introduce the spectral parameter $\lambda$ as
mentioned above we end with the two matrices:
\begin{equation}
  \label{eq:volterraLax}
  \begin{array}{rcl}
    L^{\mathrm{v}}(\lambda) &=&
    \quadmatrix{1}{\lambda\q_k}{-1}{0}\\[0.4cm]
    V^{\mathrm{v}}(\lambda) &=&
    \quadmatrix{1+\lambda\q_{k+1}}{\lambda\q_{k+1}}{-1}{\lambda\q_k}
  \end{array}
\end{equation}
with the compatibility condition 
$\lambda\dot L^{\mathrm{v}}_k(\lambda) = V^{\mathrm{v}}_k(\lambda) L^{\mathrm{v}}_k(\lambda) -
L^{\mathrm{v}}_k(\lambda) V^{\mathrm{v}}_{k-1}(\lambda)$.
This is up to the change $\lambda\to\lambda^{-2}$ and a gauge
transformation with 
$$\tilde E =
\quadmatrix{\lambda^{1/2}}{0}{0}{\lambda^{-1/2}}
$$
 the known form of
the Volterra Lax-pair \cite{SU02}.

%%%
%%% Euclidean reduction
%%%
\subsection{Euclidean reduction}
\label{subsec:dc:discreteEuclidean}

Let $c$ be a discrete curve in $\C$ and set $S_k :=
c_{k+1} - c_k$. If $\vert S_k\vert=1$ holds we call $c$ arc length
parametrized. 

Remember, that a flow on $c$ can be described via
lemma~\ref{thm:dc:cp1Flow}
For an arc length parameterized curve $c$ the {\em curvature} $\kappa$
is defined as follows:
\begin{equation}
  \label{eq:dc:curvature}
  \kappa_k = 2\tan\frac{\angle(S_{k-1},S_k)}2\thep
\end{equation}
$\kappa$ can be computed in the following way:
\begin{equation}\label{eq:dc:curvature2}
  \kappa_k = 2 i \frac{1 - \frac{S_k}{S_{k-1}}}{1 +
  \frac{S_k}{S_{k-1}}}
\end{equation}

In the arc length parametrized  case we can write
\[ 
M^h(c_{k+1}-c_k,c_k-c_{k-1}) = 
 2\frac{S_{k-1} S_k}{S_{k-1}+S_k} 
= \frac{S_{k-1} + S_k}{1 + \<S_{k-1},S_k>}
\]
since 
$$
\<2\frac{S_{k-1} S_k}{S_{k-1}+S_k},S_{k-1}> 
= \Re(2\frac{S_{k-1} S_k}{S_{k-1}+S_k}\bar
S_{k-1}) = 1 = \<\frac{S_{k-1} + S_k}{1 + \<S_{k-1},S_k>},S_{k-1}>
$$
and the same for the scalar product with $S_k$.

Let us compute how the discrete curvature evolves:
Write
\def\MU{\mu}
\begin{equation}
  \label{eq:dc:dotS}
  \dot S_k = i \MU_k S_k
\end{equation}
and since then
$\frac{d}{dt}\frac{S_k}{S_{k-1}} =
(\MU_k-\MU_{k-1}) i\frac{S_k}{S_{k-1}}$ we get
\begin{equation}
  \dot\kappa = 2 i
  \left(
  \frac{d}{dt}\frac{S_k}{S_{k-1}}\right) \frac{-2}{(1+
  \frac{S_k}{S_{k-1}} )^2} = (\mu_k-\mu_{k-1}) (1 + \frac{\kappa_k^2}4)
\end{equation}
On the other hand using equations (\ref{eq:dc:dotS}) and
(\ref{eq:dc:curvature2})  one can calculate 
$\MU_k$ to be
\begin{equation}
  \label{eq:dc:MU}
  \MU_k = \frac12
  \left(
    \beta_{k+1}\kappa_{k+1} - \beta_k\kappa_k
  \right) 
  -i
  \left(
     \beta_{k+1} + \beta_k
  \right)
\end{equation}
In the case $\beta = 2$ this implies for the evolution of the discrete
curvature $\kappa$ 
\begin{equation}
  \label{eq:dc:kappaTangential}
  \frac{\dot\kappa_k}{1+\frac{\kappa_k^2}4} = \kappa_{k+1} - \kappa_{k-1}
\end{equation}
and as the flow on the discrete curve we get the well known
tangential flow \cite{DS99,BS98}:
\[\dot c = \frac{S_- + S}{1 + \<S_-,S>}\thep\]

Now let us rewrite $\q$ to get an interpretation for the choice of
$\beta$ that gives the second Volterra flow 
$\beta_k = (\q_{k-1}+\q_{k}+1)/2$:

\[ 
\begin{array}{rcl}
\q_k &=& \frac{S_{k-1} S_{k+1}}{(S_{k-1} + S_k)(S_k+S_{k+1})} 
= ((1+\frac{S_k}{S_{k-1}})(\frac{S_k}{S_{k+1}}+1))^{-1}\\[0.4cm]
&=& ((1+\frac{2i-\kappa_k}{2i+\kappa_k})(1+\frac{2i+\kappa_{k+1}}{2i-\kappa_{k+1}}))^{-1}
= -\frac1{16}(2i+\kappa_k)(2i-\kappa_{k+1})\\[0.4cm]
&=& \frac1{16}(2i(\kappa_{k+1} -\kappa_k) + (\kappa_k \kappa_{k+1}) + 4)\thep
\end{array}
\]
With this on hand we can calculate
\[
\MU_k = \frac{11}{32}(\kappa_{k+1}+\kappa_k) + \frac{1}{64}
\left(
  (\frac{\kappa_{k+1}^2}{4}+1)(\kappa_{k+2} + \kappa_k) +
  (\frac{\kappa_k^2}{4} + 1) (\kappa_{k+1} + \kappa_{k-1})
\right)
\]
which leaves us with
\begin{equation}
  \label{eq:dc:discreteMKdV}
  \begin{array}{rcl}
    \frac{\displaystyle\dot\kappa_k}{1+\frac{\kappa_k^2}{4}} &=& 
    \frac{11}{32}(\kappa_{k+1}-\kappa_{k-1}) +\\[0.2cm]
    &&\frac{1}{64} \left(
      (\frac{\kappa_{k+1}^2}{4}+1)(\kappa_{k+2} + \kappa_k) -
      (\frac{\kappa_{k-1}^2}{4} + 1) (\kappa_{k} + \kappa_{k-2})
    \right)
  \end{array}
\end{equation}
for the evolution of the discrete curvature $\kappa$.  This is---up to
a tangential flow part which can be removed by adjusting the constant
term in the choice of $\beta$---a discretization of the mKdV equation:
\begin{equation}
  \label{eq:dc:smoothMKdVeq}
  \dot \kappa = \kappa''' +\frac32\kappa^2\kappa'\thep
\end{equation}

Therefore we will call the flow that comes from the second Volterra
flow discrete {\em mKdV flow}:
\begin{equation}
  \label{eq:dc:euclidianVL2}
  \dot c = \frac1{32}((\kappa_k \kappa_{k-1} + \kappa_{k+1} \kappa_{k}) +
  2 i (\kappa_{k+1} - \kappa_{k-1})+ 1)M^h(S_{k-1},S_k)
\end{equation}

\begin{lemma}
  The discrete tangential flow and the discrete mKdV flow both
  preserve the discrete arc length parameterization.
\end{lemma}
\begin{proof}
  We calculate $\<S_k,\dot S_k>$ for a general flow:
  \[
  \begin{array}{rcl}
    \<S_k,\dot S_k> &=& \Re(\bar S_k (\beta_{k+1}M^h(S_{k},S_{k+1}) - \beta_kM^h(S_{k-1},S_k)))\\
    &=&2 \Re\left(\frac{\beta_{k+1}}{1+\frac{S_k}{S_{k+1}}} - \frac{\beta_k}{1+ 
        \frac{S_k}{S_{k-1}}}\right)=2\Re(\beta_{k+1}(1+i\kappa_{k+1})
    -\beta_k(1-i\kappa_k))\thep 
\end{array}
  \]
  So the condition for a flow of the form $\dot c = \beta M^h(S_{k-1},S_k)$ to
  preserve the discrete arc length is
  \begin{equation}
    \label{eq:dc:euklArcKond}
    \Re(\beta_{k+1}-\beta_k) = \Im(\kappa_{k+1}\beta_{k+1} + \kappa_k \beta_k)\thep
  \end{equation}
  for the tangential flow this clearly holds. In the case of the mKdV
  flow it is an easy exercise to show equation (\ref{eq:dc:euklArcKond}).
\end{proof}

In section~\ref{sec:dc:elastica} we will see that the discrete mKdV
flow is connected to so called discrete elastic curves.

%%%
%%% Discrete Flows
%%%
\subsection{Discrete flows}
\label{subsec:dc:discreteFlows}

As in the previous section let $\g$ be the lift of a immersed discrete
curve in $\CP^1$ into $\C^2$ satisfying the normalization
(\ref{eq:dc:arcLengthConstraint}).
\begin{lemma}
  Given an initial $\tilde\g_0$ and a complex parameter $\mu$ there is
  an unique map $\tilde\g:\Z\to\C^2$ satisfying normalization
  (\ref{eq:dc:arcLengthConstraint}) and 
  \begin{equation}
    \label{eq:dc:crEvol}
    \mu = \crt(\g_k,\g_{k+1},\tilde\g_{k+1},\tilde\g_k)\thep
  \end{equation}
  We will call $\tilde\g$ a B\"acklund transform of $\g$.
\end{lemma}
  \begin{proof}
    Solving equation (\ref{eq:dc:crEvol}) for $\tilde\g_{k+1}$ gives that
    $\tilde\g_{k+1}$ is a M\"obius transform of $\tilde\g_k$.
  \end{proof}

\begin{lemma}
  If $\tilde\g$ is a B\"acklund transform of $\g$ with parameter $\mu$
  then
  \begin{equation}
    \label{eq:dc:qEvol}
    \begin{array}{rcl}
      \displaystyle  \tilde \q_k &=&\displaystyle \q_k\frac{s_k}{s_{k+1}},\\
    \displaystyle  (1-\mu) \q_k &=&\displaystyle \frac{s_{k+1}}{(1 -s_k)(s_{k+1} -1)}
    \end{array}
  \end{equation}
  with $s_k = \crt(\g_{k-1},\tilde\g_k,\g_{k+1},\g_k)$.
\end{lemma}
\begin{proof}
  Due to the properties of the cross-ratio (a useful table of the identities  can be found in \cite{HHP99}) we have
  \begin{eqnarray}
1 - \mu &=& \crt(\g_k,\tilde\g_{k+1},\g_{k+1},\tilde\g_k) = 
  \frac{\det(\g_k,\tilde\g_{k+1})}{\det(\tilde\g_{k+1},\g_{k+1})}
  \frac{\det(\g_{k+1},\tilde\g_k)}{\det(\tilde\g_k,\g_k)}\nonumber\\
  \q_k &=&\frac{\det(\g_{k-1},\g_k)}{\det(\g_k,\g_{k+2})}
  \frac{\det(\g_{k+2},\g_{k+1})}{\det(\g_{k+1},\g_{k-1})}\nonumber\\
  \frac1{1-s_k} &=& \crt(\g_{k-1},\g_k,\tilde\g_k,\g_{k+1}) =
  \frac{\det(\g_{k-1},\g_k)}{\det(\g_k,\tilde\g_k)}
  \frac{\det(\tilde\g_k,\g_{k+1})}{\det(\g_{k+1},\g_{k-1})}\nonumber\\
  \frac{s_{k+1}}{s_{k+1} -1} &=& \crt(\g_k,\tilde\g_{k+1},\g_{k+1},\g_{k+2}) =
  \frac{\det(\g_k,\tilde\g_{k+1})}{\det(\tilde\g_{k+1},\g_{k+1})}
  \frac{\det(\g_{k+1},\g_{k+2})}{\det(\g_{k+2},\g_k)}\thep\nonumber
  \end{eqnarray}
  Multiplying the first two and the second two equations proves the
  second statement.  If we set $\tilde s_k =
  \crt(\tilde\g_{k-1},\g_k,\tilde\g_{k+1},\tilde\g_k)$ we see that
  $\frac{s_k}{\tilde s_k} = 1$ and therefore
  $$(1 - \mu)\tilde \q_k = \frac{\tilde s_{k+1}}{(1 -\tilde s_k)(\tilde s_{k+1}
    -1)} = \frac{\frac1{s_{k+1}}}{(1 -\frac1{s_k})(\frac1{s_{k+1}} -1)} =
  \frac{s_k}{(1 -s_k)(s_{k+1} -1)}$$
  which proves the first statement.
\end{proof}
If $c$ is a periodic curve with period $N$, we can ask for $\tilde c$
to be periodic too. Since the map sending $c_0$ to $c_N$ is a M\"obius
transformation it has at least one but in general two fix-points. These
special choices of initial points give two B\"acklund transforms that
can be viewed as past and future in a discrete time evolution.

We will now show, that this B\"acklund transformation can serve as a
discretization of the tangential flow since the evolution on the $\q$'s
are a discrete version of the Volterra model.

The discretization of the Volterra model first appeared in
Tsujimoto, e.\ al.\ 1993.%\cite{TS93}
We will refer to the version stated in
\cite{SU02}.  There it is given in the form
\begin{eqnarray}
  \label{eq:dc:ddVolterra}
  \tilde \alpha_k = \alpha_k \frac{\beta_{k+1}}{\beta_k}\\
  \beta_k - h \alpha_k = \frac{\beta_{k-1}}{\beta_{k-1} - h \alpha_{k-1}}
\end{eqnarray}
with $h$ being the discretization constant.
\begin{theorem}
  Let $\tilde \q$ be a B\"acklund transform of $\q$ with parameter
  $\mu$.  The map sending $\q_k$ to $\tilde \q_{k+1}$ is the discrete time
  Volterra model (\ref{eq:dc:ddVolterra}) with $\alpha_k = \q_k$, $\tilde
  \alpha_k = \tilde \q_{k+1}$, $\frac{\beta_k}{h} = \frac{\q_k}{s_{k+1}}$
  and $h = \mu - 1$.
\end{theorem}
\begin{proof}
  With the settings from the theorem we have
  \[\tilde \alpha_k =\tilde \q_{k+1} = \q_{k+1}\frac{s_{k+1}}{s_{k+2}} =
  \q_k \frac{ \q_{k+1} s_{k+1}}{s_{k+2} \q_k} = \alpha_k \frac{\beta_{k+1}}{\beta_k}\]
  and on the other hand
  \[\beta_k - h \q_k = (\mu -1)\q_k(\frac1{s_{k+1}} - 1) = \frac1{1 - s_k}\]
  and \[ \frac{\beta_{k-1}}{\beta_{k-1} - h \q_{k-1}} = \frac{\frac{h \q_{k-1}}{s_k}}{h
    \q_{k-1}(\frac1{s_k} -1)} = \frac1{1 - s_k}\thep\]
This proves the theorem.
\end{proof}

The continued B\"acklund transformations give rise to maps
$\g:\Z^2\to\CP^1$ that can be viewed as discrete conformal
maps---especially in the case when $\mu$ is real negative (which is
quite far from the tangential flow, that is approximated with
$\mu\approx 1$) \cite{BP96,BP99,HMNP01}.

On the other hand in case or real $\mu$ the transformation is not
restricted to the plane: Four points with real cross-ratio always lie
on a circle. Thus the map that sends $\tilde \g_k$ to $\tilde \g_{k+1}$ is
well defined in any dimension. Maps from $\Z^2$ to $\R^3$ with
cross-ratio -1 for all elementary quadrilaterals\footnote{More general
  one can demand $\crt = \frac{\alpha_n}{\beta_m}$---see
  Chapter~\cite{BP96}.} serve as discretization of isothermic
surfaces and have been investigated in \cite{BP96}. 

\bigskip

%%%
%%% Flows and the Toda lattice
%%%
\section{Flows on discrete curves and the Toda lattice}
\label{sec:dc:todaLattice}
\begin{Definition}
Let $\lambda \in \C$ be arbitrary. Define
\beqna
p_k &:=& \frac{1}{g_k g_{k-1}}u_k -\lambda \label{def:dc:p}\\
e^{\frac{-q_{k+1}+q_k}{2}} &:=& g_k
\label{def:dc:q}\eeqna
\end{Definition}
Clearly the above definitions are not unique. $p_k$ and $g_k^{-2}$
will be identified with the Flaschka-Manakov variables of the Toda
lattice hierarchy \cite{Flasch74a,Flasch74b,MAN74}. $\lambda$ will be
the corresponding spectral 
parameter. 

With the above definitions at hand we are now able to state the
following correspondence with the Toda lattice hierarchy. 
\begin{Theorem}
Denote
$$
V_k := \left( \ba{cc}v_k^{11}&v_k^{12}\\v_k^{21}&v_k^{22}\ea\right).
$$
Define
\beqn
\alpha_k := v_k^{11}+\frac{v_k^{12}u_k}{g_{k-1}+g_k}
\hspace{5mm}
\beta_k := -\frac{v_k^{12} g_{k-1}u_k}{g_{k-1}+g_k}.
\label{def:dc:alphabet}
\eeqn
By (\ref{eq:dc:generalFlow}) and with definition
(\ref{def:dc:alphabet}),
 $\alpha_k$ and 
$\beta_k$ define a certain flow on discrete curves in
$\C^2$.
Let $V_k$ be a Lax representation matrix corresponding to the n-th Toda
flow in the notations as in \cite{FT86}. Then the Lax matrices $V_k$ of the
discrete curve flow in (\ref{eq:dc:vk}) together with the definitions 
(\ref{def:dc:alphabet}) are identical to the above
$V_k.$ Hence the compatibility equation (\ref{eq:dc:comp}) for a flow
on discrete curves in $\C^2$ corresponding to the definitions (\ref{def:dc:alphabet}), is the compatibility equation
of the n-th Toda flow. 
\label{theor:dc:major}
\end{Theorem}
\begin{proof}
Setting
\beqna
-g^{-1}_{k-1}u_k^{-1}(g_{k-1}+g_k)\beta_k &\stackrel{!}{=}& v_k^{12}\\
g^{-1}_{k-1}u_{k-1}^{-1}(g_{k-2}+g_{k-1})\beta_{k-1} &\stackrel{!}{=}& v_k^{21}
\eeqna
results in the constraint
$$
v_k^{12}\stackrel{!}{=}-v_{k+1}^{21}\frac{g_k}{g_{k-1}}
$$
But this constraint is just the 22-component of the compatibility
equation (\ref{eq:dc:comp}) for general $V_k$ and the Toda $L_k$, hence
it is satisfied by all $V_k$ of the Toda hierarchy. Hence the
variables $\beta_k$ are well defined.
Likewise the second constraint obtained by setting
\beqna
\alpha_k + \frac{1}{g_{k-1}}\beta_{k}  &\stackrel{!}{=}&  v_k^{11}\\
\alpha_{k-1} - \frac{1}{g_{k-1}}\beta_{k-1}  &\stackrel{!}{=}&
v_k^{22}
\eeqna
together with the 22-component gives the 21 component of
(\ref{eq:dc:comp}). Hence the variables $\alpha_k$ are well defined.
The 11- and 12-component are giving the Toda field equations.
\end{proof}
\begin{Remark}
The above matrices 
can be regauged into the matrices 
$\hat L_k=\Omega_{k+1}^{}L_k \Omega_k^{-1}$ and $\hat V_k=
\Omega_k^{}V_k \Omega_k^{-1} +\dot\Omega_k^{} \Omega_k^{-1} $ with
$$ \Omega_k = \left( \ba{cc} e^{\frac{q_{k}}{2}}&0\\
0&-e^{\frac{q_{k-1}}{2}}\ea \right).
$$
$\hat L_k$ and  $\hat V_k$ are then the 2 by 2 matrix representation
of the usual Flaschka-Manakov matrices  \cite{SU02}.
One has $tr \; V_k = tr \; \hat V_k - \frac{1}{2}(\dot q_k + \dot
q_{k-1})$
On the other hand by the definition of the matrices $V_k$ in 
proposition \ref{def:dc:gu}
one has
\beqn
\frac{\dot g_k}{g_k}= \frac{1}{2}(\dot q_k - \dot q_{k+1})=\tr \; V_{k+1}.
\label{eq:dc:trace}\eeqn
hence
\beqn
\dot q_k =  tr \; \hat V_{k+1}.
\label{eq:dc:qdotandtrace}
\eeqn
\label{rem:dc:qs}\end{Remark}

\section{The first three Toda lattice hierachy flows and reductions of
them}\label{sec:dc:threeFlowsRed} 
\subsection{The first three Toda flows}
\label{subsec:dc:threeFlows} 
The  flow directions for a discrete curve $\g$ is
given by a specific choice of the variables $\alpha_k$, $\beta_k$ 
(compare with (\ref{eq:dc:generalFlow})). In the following
we will choose $\alpha_k$, $\beta_k$ in such a way that the corresponding
evolution for the determinants $g_k$, and $u_k$ (\ref{def:dc:gu})
is the evolution of the canonical variables of the toda lattice hierachy.
In this section we will look at the flow directions
given by  the ``first'' three flows of the Toda lattice hierachy.
\begin{Proposition}[First flow]
Using (\ref{eq:dc:generalFlow})
define the following flow directions for $\g$:
\beqna
\alpha_k^{TL1} &=& - \frac{1}{2}(p_k +
\lambda)\frac{g_{k-1}-g_k}{g_{k-1}+g_k}=
- \frac{u_k}{2g_{k-1}g_k}\frac{g_{k-1}-g_k}{g_{k-1}+g_k}\\
\beta_k^{TL1} &=& - \frac{g_{k-1}g_k(p_k + \lambda)}{g_{k-1}+g_k}= -
\frac{u_k}{g_{k-1}+g_k }.
\label{def:dc:albetTL1}\eeqna
The induced flow on the determinants $g_k$ and $u_k$ 
(\ref{def:dc:gu}) is with definitions 
(\ref{def:dc:p}), (\ref{def:dc:q}) given by the 
first Toda lattice flow.
\label{prop:dc:firstflow}
\end{Proposition}
\begin{proof}In accordance with theorem \ref{theor:dc:major} 
we obtain the following compatibility matrices:
\beqna
L_k&=&\left(\ba{cc} g_k(p_k + \lambda)&-\frac{g_k}{g_{k-1}}\\
1&0\ea\right)
\\
V_k&=&\left(\ba{cc}
-\frac{1}{2}(p_k+\lambda) & g_{k-1}^{-1}
              \\
-  g_{k-1}^{-1} & \frac{1}{2}(p_{k-1}+\lambda)\ea\right)
\eeqna
These matrices are wellknown \cite{FT86}  and give the Toda lattice
equations, which are called the first flow of the Toda lattice hierachy:
\beqna
\frac{\dot g_k}{g_k}&=& tr \; V_{k+1} = \frac{1}{2} (p_{k}-p_{k+1})\\
\dot p_k&=& g_k^{-2}-g_{k-1}^{-2}.
\eeqna
\end{proof}
\begin{Remark}Due to  remark \ref{rem:dc:qs} the
variables $\dot q_k$  evolve with $\dot q_k=p_k+\lambda.$
Hence
$$
{\ddot{q}}_k = g_k^{-2}-g_{k-1}^{-2}
$$
which is the wellknown Toda lattice equation.
\end{Remark}
Analogously the flows for the next two higher flows 
 \cite{SU02} can be determined:
\begin{Proposition}[second flow]
Define
\beqna
\alpha_k^{TL2} &=& - \frac{1}{2}(p_k^2 -
\lambda^2)\frac{g_{k-1}-g_k}{g_{k-1}+g_k}
-\frac{1}{2}(g_{k}^{-2}-2g_{k-1}^{-2}+g_{k-2}^{-2})\\
\beta_k^{TL2} &=& - \frac{g_{k-1}g_k(p_k^2 - \lambda^2)}{g_{k-1}+g_k}.
\label{def:dc:albetTL2}\eeqna
The induced flow on  $g_k$ and $p_k$ 
(\ref{def:dc:gu}) is with definitions (\ref{def:dc:p}), (\ref{def:dc:q})
given by the second flow of the  Toda lattice hierachy \cite{SU02}:
\beqna
\frac{\dot{g}_k}{g_k}&=&-\frac{1}{2}(p^2_{k+1}-p^2_k+g_{k+1}^{-2}-g_{k-1}^2)
\label{eq:dc:evgTL2}\\
\dot p_k&=& g_k^{-2}(p_{k+1}+p_k)-g_{k-1}^{-2}(p_k+p_{k-1})
\eeqna
\end{Proposition}
\begin{Remark}$\dot q_k = g_k^{-2}+g_{k-1}^{-2}+p_k^2-\lambda^2.$
\end{Remark}
\begin{Proposition}[third flow]
Define
\beqna
\alpha_k^{TL3} &=& - \frac{1}{2}\frac{g_{k-1}-g_k}{g_{k-1}+g_k}
(p_k^3+\lambda^3-2(p_k+\lambda)g_k^{-1}g_{k-1}^{-1})\\
&&-\frac{1}{2}g_{k}^{-2}(2 p_{k}+p_{k+1})  +
\frac{1}{2}g_{k-1}^{-2}(p_{k-1}+2p_k)\\
\beta_k^{TL3} &=& - \frac{g_{k-1}g_k(p_k + \lambda)}{g_{k-1}+g_k}
(g_{k-1}^{-2}+g_{k}^{-2}+p_k^{2}-\lambda p_k +\lambda^2).
\label{def:dc:albetTL3}\eeqna
The induced flow on $g_k$ and $p_k$ 
(\ref{def:dc:gu}) is with definitions (\ref{def:dc:p}), (\ref{def:dc:q})
given by the third flow of the  Toda lattice hierachy  \cite{SU02}:
\beqna
\frac{\dot{g}_k}{g_k}&=&
-\frac{1}{2}[(p_{k+1}^3+p_{k+2}g_{k+1}^{-2}+2p_{k+1}g_{k+1}^{-2}
+p_{k+1}g_{k}^{-2})\\
&& -(p_k^3+p_{k-1}g_{k-1}^{-2}+p_{k}g_{k}^{-2}+2p_{k}g_{k-1}^{-2})]\nonumber\\
\dot p_k&=&g_k^{-2}( p_{k+1}^2 +p_k^2+ p_{k+1}p_k
+g_{k+1}^{-2}+g_{k}^{-2})- \nonumber\\
&&g_{k-1}^{-2}(p_k^2+p_{k-1}^2+p_k p_{k-1}+g_{k-1}^{-2}+g_{k-2}^{-2})
\eeqna
\end{Proposition}
\begin{Remark}$\dot q_k = g_k^{-2}(2 p_k +p_{k+1})
+g_{k-1}^{-2}( p_{k-1}+2 p_k ) + p_k^3 +\lambda^3.$
\end{Remark}
%%%
%%% Geometric interpretations
%%%
\subsection{Towards geometrical interpretations of the Toda lattice}
\label{sec:dc:invariant}
In the following we look at the reduction $p_k=0$ and determine 
the geometric shape of curves, which are
invariant under the trivial Toda and the discrete KdV flow.

\subsubsection{Reduction $p_k=0$ and invariant curves under the 
trivial Toda flow}
\label{subsec:dc:ellipse} 

In the previous section one reduced the discrete
curves in $\C^2$ to socalled arclength para\-metrized curves, by setting
 $g_k =1 $ for all $k.$ It was conjectured
that for that reduction there exist certain flow directions (\ref{def:dc:betaIt})
which let the cross-ratios evolve according to flows in the 
Volterra hierachy. 

On the other hand it is a fact  
\cite{SU02} that the Toda lattice hierachy reduces to the 
Volterra hierachy for flows with an even enumeration number
if $p_k=0$. Indeed, looking at the above equations
one sees easily that if $p_k=0$ then  the evolution of the $g_k^{-2}$ 
in equation
\ref{eq:dc:evgTL2} is given by the Volterra equation.
In the same manner the third flow admits no reduction but the fourth
flow would  again give an equation for the $g_k^{-2}$ which can be identified 
with a discrete version of the KdV equation \cite{SU02}.   
We observe that for the  case  $p_k=0$ the cross-ratio is given
by 
$$
\q_k =\lambda g_k^{-2}.
$$
Hence as a direct consequence
\begin{Proposition}
If $\lambda=1$, $p_k=0$ then the flow of the  cross-ratios $\q_k$,
as given by the flow of the  2n'th Toda lattice hierachy  
via theorem \ref{theor:dc:major} satisfies the equations
of the n'th Volterra hierachy.
\end{Proposition}

As can be seen by looking at the above Toda lattice hierachy
equations, the reduction $g_k=const$ is generally not compatible with
the equations. It is compatible if $p_k=const$ (and especially
$p_k=0$, which can allways be achieved by a change of $\lambda$),
which gives the trivial evolution $\dot g_k =0, \dot p_k =0$ for all
flows of the Toda hierachy. Hence the two reductions $\dot g_k =0$ and
$p_k=0$ with their corresponding Volterra hierachy flows, seem to
be in some kind of duality, where only the trivial reduction $g_k
=const \wedge p_k=const$ $\leftrightarrow$ $u_k=const$ $\wedge$
$g_k=const$ seems to lie in the intersection of the two pictures. We
were interested in what kind of curves belong to this most trivial
solution of the Toda flows.
The following proposition shows that if the variables
$u_k=u$ and $g_k=g$ are constants, then $\g$ lies on a quadric,
or in other words in this case $\g$ defines a {\bf discrete quadric}
in $\C^2.$
\begin{Proposition}
a) Let $\g_0$, $\g_1$ be fixed initial conditions for
a discrete curve $\g:\Z\to\C^2$  with
\beqna
\<M \g_0,\g_0> &=& 1\\
\<M \g_1,\g_1> &=& 1\\
\<M \g_0,\g_1> &=& \frac{u}{2g},\\
\label{def:dc:initial1}\eeqna
where $\<\g,\g> = x^2+y^2$ is the naive complexification of the real scalar product to $\C^2$ and
$M$ is a symmetric 2 by 2 matrix; $g = \det(\g_0,\g_1) \neq 0$ and
$u$ is a free parameter. Then $M$ is uniquely defined by
conditions \ref{def:dc:initial1} and $\det M = \frac{1}{g^2}(1-\frac{u^2}{g^2})$. 

b)Define a discrete curve with the above initial conditions
recursively by:
$$
\g_{k+2}:= \frac{u}{g}\g_{k+1} - \g_k
$$
then 
\beqna
\<M \g_k,\g_{k+1}>& =&\frac{u}{2g}  \\
\<M \g_k,\g_k> &=& 1 \hspace{5mm} \mbox{for all} \hspace{5mm} k \in \Z.
\eeqna
\label{prop:dc:quadric}\end{Proposition}
\begin{proof}

\noindent a) After a lengthy calculation M can be derived from
the conditions (\ref{def:dc:initial1}) as
\beqn
M=
\frac{1}{g^2}\left( \ba{cc}y_0^2+y_1^2-2y_0 y_1 \frac{u}{2g}&
(x_0 y_1 + y_0 x_1)  \frac{u}{2g}-(x_0 y_0+ x_1 y_1)\\
(x_0 y_1 + y_0 x_1)  \frac{u}{2g}-(x_0 y_0+ x_1 y_1)&
x_0^2+x_1^2-2x_0 x_1 \frac{u}{2g} \ea\right)
\eeqn

\noindent b) $\<M \g_2,\g_1>=\frac{u}{g}\<M \g_1,\g_1> - \<M\g_0,\g_1>=
\frac{u}{g}\cdot 1- \frac{u}{2g}$, hence by induction 
 $\<M \g_k,\g_{k+1}>=\frac{u}{2g}$ for all $k \in \Z$.
From this it follows that 
$$ \<M  \g_2,\g_2>=\frac{u^2}{g^2}\<M \g_1,\g_1>+\<M \g_0,\g_0> 
-2\frac{u}{g}\<M \g_1,\g_0>=1 
$$
and hence by induction
$\<M \g_k,\g_k> =1$ for all $k \in \Z$.
\end{proof}
Looking at the variables $\alpha_k$, $\beta_k$ corresponding to the
first three Toda flows (e.g.(\ref{def:dc:albetTL1})) the case 
$p_k=const \wedge g_k = const$ gives the following evolution
on the curves:
\beqn
\dot \g_k = const(\g_{k+1}-\g_{k-1})
\label{eq:dc:tangent}\eeqn
where the constant $const$ varies corresponding to the 
considered Toda flow. 
\begin{Proposition}
A flow direction as in \ref{eq:dc:tangent} along a discrete quadric
$\g$
(as defined in proposition \ref{prop:dc:quadric}) is tangent to
the  (smooth) quadric on which the points of $\g$ lie.
\end{Proposition}
\begin{proof}
We have to show that if 
$\dot \g_k=\frac{d}{dt} \g_k = \rho(\g_{k+1}-\g_{k-1})$
then $$\frac{d}{dt}|_0 <M \g_k(t),\g_{k}(t)>=0 \hspace{0.5cm} \mbox{and}  
\hspace{0.5cm}\frac{d}{dt}|_0 \<M \g_k(t),\g_{k+1}(t)>= 0.$$ Now
 \beqnao
\frac{d}{dt}|_0 \<M \g_k(t),\g_{k}(t)>&=& 
\<M \dot \g_k,\g_{k}> +\<M \g_k,\dot \g_{k}>\\
&=& 2 \rho(\<M \g_{k+1},\g_{k}> - \<M \g_{k-1},\g_k>)\\
&=& 2 \rho(\frac{u}{2g}-\frac{u}{2g})
\eeqnao
Using lemma \ref{lem:dc:Gam1} the proof for the second assertion works 
analogously.
\end{proof}
By the above, a class of curves (namely quadrics) was sorted out
by looking at trivial solutions to the Toda lattice equations.
The movement of these curves under the corresponding Toda flow  (tangent
to the quadric) was very geometrical and simple.

In the next section a similar construction will be done. Here
it will be shown in the euclidian reduction
(\ref{subsec:dc:discreteEuclidean}),
 that discrete curves
which evolve under the discrete mKdV flow simply by a translation,
define socalled discrete elastic curves. 

\subsubsection{Invariant curves under the discrete mKdV flow: Discrete generalized elastic curves}
\label{sec:dc:elastica}

Let $c$ be a discrete arc length parametrized curve in $\C$ as
discussed in section~\ref{subsec:dc:discreteEuclidean}.

\begin{definition}
  A discrete regular arc length parametrized curve $c:\{0,\ldots,N\}$
  $\to\C$ is called \emph{planar elastic curve} if it is an critical
  point to the functional
  \begin{equation}
    \label{eq:dc:discreteElasticFunctional}
    \sum_{j=1}^{N-1} \log(1+\frac{\kappa_j}{4})
  \end{equation}
  The admissible variations preserve the arc length, $c_N-c_0$ and
  the tangents at the end points. 
\end{definition}
It can be shown \cite{BHS03} that the curvature of a discrete elastic
curve obeys the following equation:
\begin{equation}
  \label{eq:dc:elasticCurvature}
  \kappa_{k+1} = 2a \frac{\kappa_k}{1+\frac{\kappa_k^2}{4}} -
  \kappa_{k-1}\thep
\end{equation}
for some real constant $a$.

We now want to know, which class of curves is invariant under
the discrete mKdV flow (up to
euclidean motion and some tangential flow). It will turn out
that that discrete elastic curves are a special case in
that class. 

Since discrete arc length parametrized curves are determined by their
curvature (\ref{eq:dc:curvature}) up to euclidean motion, it is
sufficient to impose the constraint that the curvature must not change
up to the changes made by the tangential flow
(\ref{eq:dc:kappaTangential}). In other words: we ask for solutions
$\kappa$ for which (\ref{eq:dc:kappaTangential}) is a multiple of
(\ref{eq:dc:discreteMKdV}):
\begin{equation}
  \label{eq:dc:kappaConstraint}
  (a-\frac{11}{32})(\kappa_{k+1} - \kappa_{k-1}) = 
  \frac{1}{64} \left(
    (\frac{\kappa_{k+1}^2}{4}+1)(\kappa_{k+2} + \kappa_k) -
    (\frac{\kappa_{k-1}^2}{4} + 1) (\kappa_{k} + \kappa_{k-2})
  \right)
\end{equation}
One can ``integrate'' this equation twice and get the following lemma:
\begin{lemma}
  The curvature of a discrete curve, that evolves up to some
  tangential flow by Euclidean motion under the mKdV flow satisfies 
  \begin{equation}
    \label{eq:dc:gElastic}
    \kappa_{k+1} = \frac{2a \kappa_k}{1 + \kappa_k^2} - \kappa_{k-1} +
    b + c_k
  \end{equation}
  for some constants $a$ and $b$ and a function $c$ with $c_{k+1} = -c_k$.
\end{lemma}
\begin{proof}
  Nothing left to show.
\end{proof}

\begin{figure}[htbp]
 \begin{center}
   \epsfxsize=0.35\hsize\epsfbox{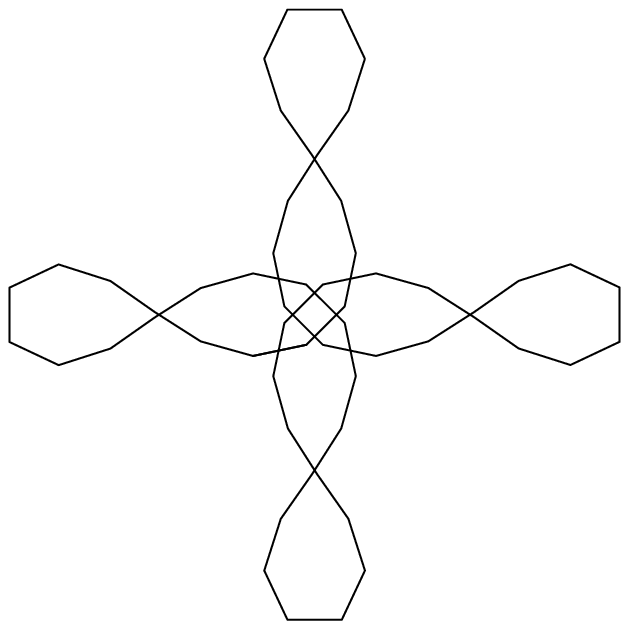}\hfil\epsfxsize=0.35\hsize\epsfbox{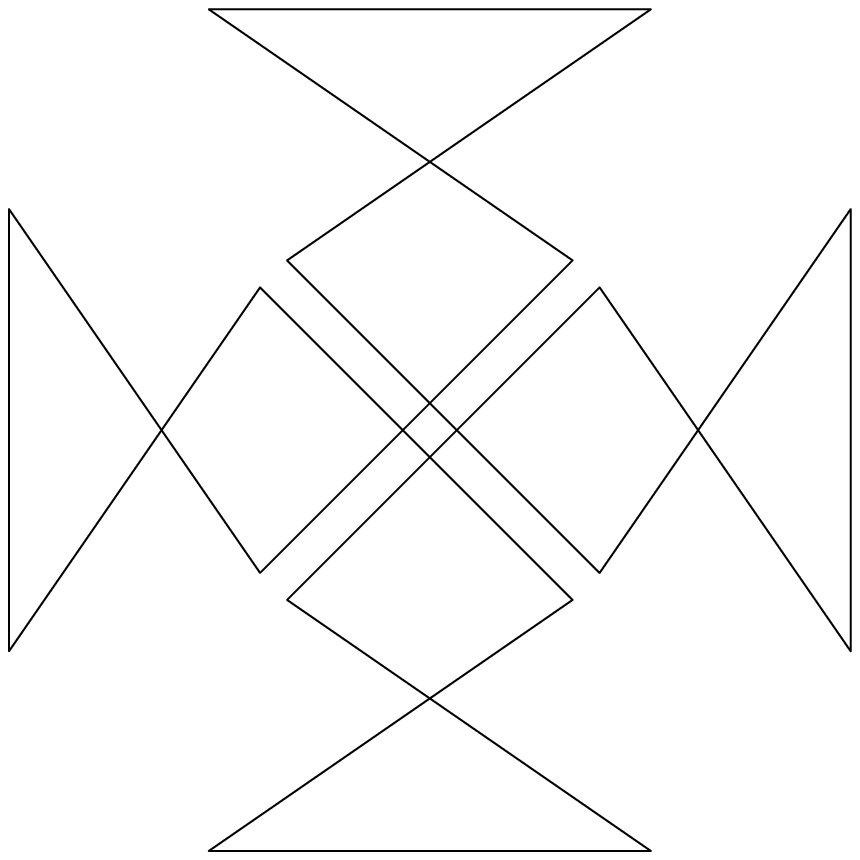}
\centerline{\epsfxsize=0.35\hsize\epsfbox{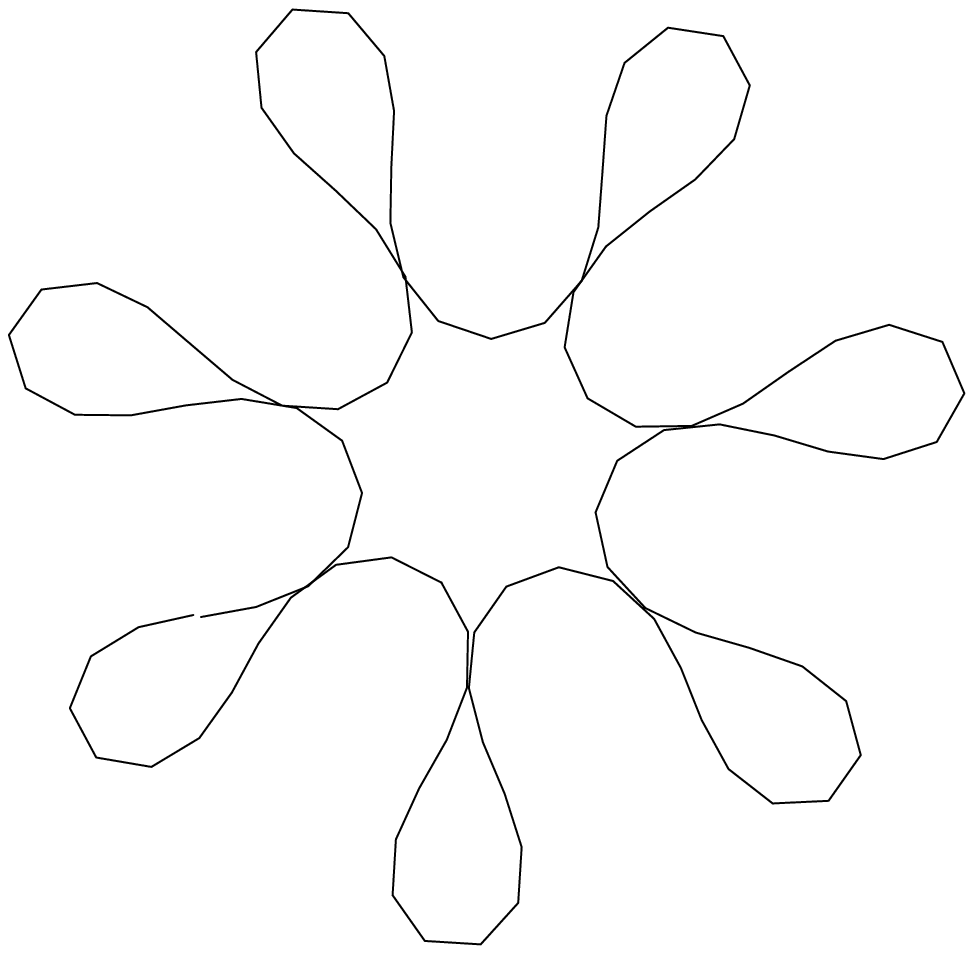}}
   \caption{Generalized elastic curves.}
   \label{fig:dc:fourSym}
 \end{center}
\end{figure}

%\begin{remark}
In the case $b = 0$ and $c\equiv 0$ this gives the equation for
planar elastic curves (\ref{eq:dc:elasticCurvature}).
%\end{remark}
Figure~\ref{fig:dc:fourSym} shows three closed discrete generalized
elastic curves.

\section{Conclusion}
In this paper we introduced a novel view onto the one dimensional
Toda lattice hierarchy in terms of flows on discrete curves. In particular 
this view allowed us to give a geometric meaning to the trivial Toda flow 
and the mKdV flow which are special flows within the Toda lattice hierarchy. 
This was achieved by classifying that class of curves, which is invariant 
under the corresponding flow. 
It would be interesting to investigate discrete curves that belong to other
Toda flows in this sense. In
fact we view
our exposition rather as a starting point for a more thorough investigation
of this subject.

If one finds a meaningful symplectic structure on the space of
curves $\gamma$, then by the novel view onto the  Toda lattice hierarchy this 
may lead to a symplectic structure for the ``vertex operators'' $F_n$.
We are currently working on that issue.

The theory of discrete curves is in close connection to the
theory of discrete surfaces. 
In addition it was already known to Darboux \cite{Dol97} that 
certain invariants 
on surfaces admitting a conjugate net
parametrization satisfy the two dimensional Toda lattice equation.
It is an  interesting question to see whether
our approach can be transferred to discrete surfaces and the
two dimensional Toda lattice hierarchy \cite{UeT84}.  

It would be interesting to see, whether our approach could also
be applied to generalized Toda Systems \cite{Ko79}. We haven't 
thought about that yet.

\section*{Acknowledgements}
We like to thank Ulrich Pinkall for deep and interesting discussions
and many inspiring ideas. We like to thank Yuri Suris and Leon Takhtajan 
for helpful hints.

%\bibliographystyle{alpha}
%\bibliography{discrete}

\end{document}